\documentclass[11pt,leqno]{article} 
\usepackage{graphics}
\newtheorem{thm}{Theorem}[section]
\newtheorem{lma}{Lemma}[section]
\newtheorem{cor}{Corollary}

\newcommand{\beqa}{\begin{eqnarray}}
\newcommand{\eeqa}{\end{eqnarray}}

\newcommand{\pf}{\noindent {\bf Proof:} $\s$ }
\newcommand{\epf}{ \hfill$\diamondsuit$ \medskip}
\newcommand{\md}{\medskip}
\newcommand{\ds}{\displaystyle}
\newcommand{\beq}{\begin{equation}}
\newcommand{\eeq}{\end{equation}}
\newcommand{\lbl}{\label}
\newcommand{\s}{\; \;}

\newcommand{\la}{\lambda}

\newcommand{\ra}{\rightarrow}

\title{A global solution curve for a class of free boundary value problems arising in plasma physics}

\author{
Philip Korman   \\ 
Department of Mathematical Sciences \\ 
University of Cincinnati \\ 
Cincinnati Ohio 45221-0025 \\
}

\date{}

\begin{document}

\maketitle
\begin{abstract} 
We   study the existence and multiplicity of solutions and the global solution curve    of the following free boundary value problem, arising in plasma physics, see H. Berestycki and H. Brezis \cite{BB}: find a function $u(x)$ and a constant $b$, satisfying
\beqa \nonumber
& \Delta u+g(x,u)=p(x) \s \mbox{in $D$} \\ \nonumber
& u \,| \, _{\partial D}=b, \s\s \int_{\partial D} \frac{\partial u}{\partial n} \, ds=0 \,. \nonumber
\eeqa
Here $D \subset R^n$, is  a bounded domain, with  smooth boundary. This problem can be seen as a PDE generalization of the periodic problem for  one-dimensional pendulum-like equations. We use continuation techniques. Our approach is suitable for numerical computations.
 \end{abstract}

\begin{flushleft}
\end{flushleft}

\begin{flushleft}
\end{flushleft}

\section{Introduction}
\setcounter{equation}{0}
\setcounter{thm}{0}
\setcounter{lma}{0}

We  study the existence and multiplicity of solutions, and the global solution structure   of the following free boundary value problem, arising in plasma physics, see H. Berestycki and H. Brezis \cite{BB}. On a  bounded domain $D \subset R^n$, with a smooth boundary, we seek to  find a function $u(x)$ and a constant $b$, satisfying
\beqa
\lbl{p1}
& \Delta u+g(x,u)=p(x) \s \mbox{in $D$} \\
& u \,| \, _{\partial D}=b, \s\s \int_{\partial D} \frac{\partial u}{\partial n} \, ds=0 \,. \nonumber
\eeqa
Here $g(x,u)$ and $p(x)$ are given functions, $n$ denotes the unit normal vector on $\partial D$, pointing outside. More recently, this problem was studied in P. Amster, P.L. De Napoli, and M.C. Mariani \cite{amster}. As explained in \cite{amster}, this problem can be seen as a PDE generalization of the periodic problem for the one-dimensional pendulum-like equation:
\beqa
\lbl{p2}
& u''+g(t,u)=p(t) \\
& u(0)=u(T)=b, \s u'(0)=u'(T) \,.  \nonumber
\eeqa
with $p(t)$ being $T$-periodic, and $g(t,u)$ is $T$-periodic in $t$. Indeed, writing the second boundary condition in (\ref{p2}) as $\int_0^T u'' \,dt=0$, we see that a similar condition for (\ref{p1}) is $\int _D \Delta u \, dx=0$, which  by the divergence theorem is equivalent to the second boundary condition in (\ref{p1}). There exists an enormous literature for the periodic problem (\ref{p2}), see e.g., the review paper of J. Mawhin \cite{M}, and \cite{D}, \cite{F}, \cite{K}, \cite{R}, \cite{R1}, \cite{T}. By constrast, not much is known for the more general problem (\ref{p1}).
\medskip

We  use continuation techniques to study the solution curves for the problem (\ref{p1}), similarly to our approach to the pendulum-like equations \cite{K2}. One can think of the problem (\ref{p1}) as being  ``at resonance", i.e., some conditions on $p(x)$ are necessary for its solvability. Indeed, decompose $p(x)=\mu _0 +\theta (x)$, with  $ \mu _0=\int_{D} p(x) \, dx$, and $\int_{D} \theta(x) \, dx=0$. Integrating the equation in (\ref{p1}), we get
\beq
\lbl{p3}
\int_{D} g(x,u(x)) \, dx=\mu _0 |D| \,,
\eeq
where $|D|$ is the volume of $D$. Assume that the following limits exist: $\ds g(x, \pm \infty)=\lim _{u \ra \pm \infty} g(x,u) $, uniformly in $x \in \bar D$, with $g(x, \pm \infty) \in L^{\infty} (D)$,  and suppose that
\beq
\lbl{p4}
g(x, - \infty)<g(x,u)<g(x, \infty) \s\s \mbox{for all $x \in D$} \,.
\eeq
Then, from (\ref{p3}) and (\ref{p4}), we obtain 
\beq
\lbl{p5}
\frac{1}{|D|} \int_{D} g(x,-\infty) \, dx<\mu _0 < \frac{1}{|D|}\int_{D} g(x,\infty) \, dx
\eeq
to be a necessary condition for existence of solutions. Similarly to the classical paper of E.M. Landesman and A.C. Lazer \cite{L}, we give conditions under which the condition (\ref{p5}) is sufficient for existence.
\md 

In order to use continuation techniques, we will embed the problem (\ref{p1}) into a family of problems
\beqa
\lbl{p6}
& \Delta u+kg(x,u)=p(x) \s \mbox{in $D$} \\
& u \,| \, _{\partial D}=b, \s\s \int_{\partial D} \frac{\partial u}{\partial n} \, ds=0 \,, \nonumber
\eeqa
with $0 \leq k \leq 1$. At $k=1$, one recovers the original problem (\ref{p1}), while at $k=0$ the equation is linear. It turns out that at $k=0$, the problem has infinitely many solutions, and one may fix a unique solution, by additionally prescribing any value of $\frac{1}{|D|} \int_{D} u(x) \, dx \equiv \xi _1$. One then performs continuation in $k$ on curves of fixed average. Namely, for any $\xi _1 \in R$, one solves the following problem: find $(u,b,\mu)$, as a function of $k$, solving
\beqa
\lbl{p7}
& \Delta u+kg(x,u)=\mu+\theta (x) \s \mbox{in $D$} \\ \nonumber
& u \,| \, _{\partial D}=b, \s\s \int_{\partial D} \frac{\partial u}{\partial n} \, ds=0 \\ \nonumber
& \frac{1}{|D|} \int_{D} u(x) \, dx = \xi _1 \,. \nonumber
\eeqa

To prove the solvability of (\ref{p7}), we set up the appropriate function spaces, and show that the corresponding differential operator is one-to-one, and onto, so that the Implicit Function Theorem applies. ``Onto" is the harder part.  Once the continuation process is completed, and one has a solution of the problem (\ref{p7}) at $k=1$, it remains to show that one can select a value of $\xi _1$, at which $\mu =\mu _0$, thus giving us a solution of the original problem (\ref{p1}). This part is also accomplished by continuation.
\medskip 

In addition to a result of Landesman-Lazer type, we obtain an existence result of D.G. de Figueiredo and W.-M.  Ni \cite{FN} type, which does not require that limits at infinity exist. We show that  $\xi _1$ is a global parameter, and then we study the curve $\mu =\mu (\xi _1)$, yielding  a multiplicity result. The continuation approach of this paper is similar to our paper on pendulum-like equations, see \cite{K2}. In addition to its conceptual clarity, this approach is suitable for  numerical computation of all solutions of (\ref{p1}). Every solution can be obtained by two continuations, first in $k$, $0 \leq k \leq 1$, and then in $\xi _1$. These solutions curves are easy to compute numerically e.g., by the predictor-corrector method, since we prove that no turning points (or other singularities) are encountered. We had performed similar numerical computations in \cite{K1}.

\section{Continuation of solutions of any fixed average}
\setcounter{equation}{0}
\setcounter{thm}{0}
\setcounter{lma}{0}
\setcounter{cor}{0}

We begin with the following linear problem
\beqa
\lbl{8}
& \Delta u=\theta (x) \s \mbox{in $D$} \\ \nonumber
& u \,| \, _{\partial D}=b, \s\s \int_{\partial D} \frac{\partial u}{\partial n} \, ds=0   \\ \nonumber
&  \frac{1}{|D|} \int _D u(x) \, dx=\xi _1 \,. \nonumber
\eeqa

\begin{lma}\label{lma:1}
Given any $\theta (x) \in L^2(D)$, with $\int _D \theta (x) \, dx=0$, and any $\xi _1 \in R$, we can find a unique pair $(u,b) \in W^{2,2}(D) \cap W_0^{1,2}(D) \times R$, solving (\ref{8}).
\end{lma}

\pf
Let $v(x) \in W^{2,2}(D) \cap W_0^{1,2}(D)$ be the solution of
\[
\Delta v=\theta (x) \s \mbox{in $D$}, \s v=0 \s \mbox{on $\partial D$}.
\]
By the divergence theorem, $ \int_{\partial D} \frac{\partial v}{\partial n} \, ds=0$. Then $u=v+b$ will give us solution of (\ref{8}), if we select the constant $b$, so that $\frac{1}{|D|} \int _D (v+b) \, dx=\xi _1$. Uniqueness follows from the fact that any harmonic in $D$ function, satisfying $\int_{\partial D} \frac{\partial u}{\partial n} \, ds=0$, is a constant.
\epf

\begin{cor}\lbl{cor:1}
Assume in addition that $\theta (x) \in L^{\infty} (D)$, and $\xi _1 =0$. Then there is a constant $c$, such that 
\[
||u||_{L ^{\infty} (\bar D) } \leq c \,.
\]
\end{cor}

\pf
Using the $W^{2,p}$ estimates and the Sobolev embedding theorem, we conclude an estimate for $||v||_{L ^{\infty} (\bar D) } $, and then for $|b|$, giving us an estimate of $u=v+b$.
\epf

The following Poincare's inequality is well known. By $c_0$ we denote the best (largest) constant. 
\begin{lma}\label{lma:0}
Assume that $u(x) \in W_0^{1,2}(D)$ is any function satisfying $\int _D u(x) \, dx=0$.
Then there is a constant $c_0$, depending only on $D$, such that
\[
\int _D |\nabla u|^2 \, dx \geq c_0 \int _D u^2 \, dx \,.
\]
\end{lma}

We consider next the following linear problem: given the function $a(x)$, find a triple $(w(x), \mu ^*,b) \in W^{2,2}(D) \times R \times R$, satisfying
\beqa
\lbl{9}
& \Delta w+a(x)w=\mu ^*, \s \mbox{in $D$} \\ \nonumber
& w \,| \, _{\partial D}=b, \s\s \int_{\partial D} \frac{\partial w}{\partial n} \, ds=0   \\ \nonumber
& \int _D w(x) \, dx=0 \,. \nonumber
\eeqa

\begin{lma}\label{lma:2}
Assume that
\[
a(x)<c_0, \s \mbox{for all $x \in D$}.
\]
Then the only solution of the problem (\ref{9}) is $w(x) \equiv 0$, $\mu ^*=0$, $b=0$.
\end{lma}

\pf
Multiply the equation in (\ref{9}) by $w$, and integrate over $D$. Since 
\[
-\int _D w \Delta w  \, dx=\int _D  |\nabla w|^2 \, dx+b \int_{\partial D} \frac{\partial w}{\partial n} \, ds=\int _D  |\nabla w|^2 \, dx \,,
\]
we have by Poincare's inequality
\[
c_0 \int _D w^2 \, dx > \int _D a(x)w \, dx =-\int _D w \Delta w  \, dx=\int _D  |\nabla w|^2 \, dx \geq c_0 \int _D w^2 \, dx \,,
\]
from which we conclude that $w(x) \equiv 0$, and then $\mu ^*=0$, and $b=0$.
\epf

\begin{thm}\lbl{thm:1}
Consider the problem
\beqa
\lbl{10}
& \Delta u+kg(x,u)=\mu+\theta (x) \s \mbox{in $D$} \\ \nonumber
& u \,| \, _{\partial D}=b, \s\s \int_{\partial D} \frac{\partial u}{\partial n} \, ds=0 \\ \nonumber
& \frac{1}{|D|} \int_{D} u(x) \, dx = \xi _1 \,, \nonumber
\eeqa
where $\theta (x) \in L^2(D)$ is a given function. The function $u(x)$ and the constants $\mu$ and $b$ are unknown. Assume that $g(x,u) \in C^1(\bar D \times R)$ satisfies
\beq
\lbl{11}
|g_u(x,u)| \leq M, \s \mbox{for all $x \in \bar D$, and $u \in R$} \,,
\eeq
and 
\beq
\lbl{12}
M < \min (c_0, \la _2) \,.
\eeq
Then for any $0 \leq k \leq 1$, and $\xi _1 \in R$, there exists a unique solution $(u(x),b,\mu) \in  W^{2,2}(D)  \times R \times R$, solving (\ref{10}).
\end{thm}

\pf
We perform continuation in $k$. When $k=0$, we take $\mu=0$, and obtain the unique solution $(u(x),0, b)$ by Lemma \ref{lma:1}.
We consider two cases.
\md 

\noindent
{\bf Case I}: $\xi _1=0$.
\md 

We shall recast the problem (\ref{10}) in the operator form, after we define the appropriate spaces. We shall denote by $W^{2,2}(D) \cap W_0^{1,2}(D)+R$ the functions of the form $v(x)+b$, with $v(x) \in W^{2,2}(D) \cap W_0^{1,2}(D)$, and $b \in R$. Define
\[
\bar H^2= \left\{ u(x) \in W^{2,2}(D) \cap W_0^{1,2}(D)+R \;\, | \; \int _D u(x) \, dx=0, \; \int_{\partial D} \frac{\partial u}{\partial n} \, ds=0\right\} \,.
\]
One checks that $\bar H^2$ is a Banach space. Consider the map
\[
F(u,\mu ,k)=\Delta u +k g(x,u)-\mu \; : \; \bar H^2 \times R  \times R \ra L^2(D) \,.
\]
Then the problem (\ref{10}) can be written in the operator form
\beq
\lbl{14}
F(u,\mu ,k)=\theta (x) \,.
\eeq
Observe that the unknown constant $b$ is now ``placed" in the definition of the space $\bar H^2$.
We wish to apply the Implicit Function Theorem to perform continuation in $k$. For that we need to show that the linearized operator
\[
F_{(u,\mu)}(u,\mu ,k)(w, \mu ^*)=\Delta w+kg_u(x,u)w-\mu ^*
\]
is one-to-one and onto.
\md 

It is one-to-one, because the only solution of the problem 
\beqa
\lbl{15}
& \Delta w+kg_u(x,u)w-\mu ^*=0, \s \mbox{in $D$} \\ \nonumber
& w \,| \, _{\partial D}=b, \s\s \int_{\partial D} \frac{\partial w}{\partial n} \, ds=0   \\ \nonumber
& \int _D w(x) \, dx=0  \nonumber
\eeqa
is $w(x) \equiv 0$, $ \mu ^*=0$, in view of Lemma \ref{lma:2}.
\md 

Turning to the ``onto" part, we need to show that given any $\theta (x) \in L^2(D)$, one can find $w \in \bar H^2$, and $\mu ^* \in R$, solving
\beqa
\lbl{16}
& \Delta w+kg_u(x,u)w=\mu ^*+\theta (x), \s \mbox{in $D$} \\ \nonumber
& w \,| \, _{\partial D}=b, \s\s \int_{\partial D} \frac{\partial w}{\partial n} \, ds=0   \\ \nonumber
& \int _D w(x) \, dx=0 \,. \nonumber
\eeqa
We consider two further cases.
\md 

\noindent
{\bf Case 1}. The operator $L \; : \; W^{2,2}(D) \cap W_0^{1,2}(D) \ra L^2(D)$, defined by
\[
Lw \equiv \Delta w+kg_u(x,u)w, \s \mbox{subject to $w=0$ on $\partial D$} 
\]
is invertible. Let $w_1 =L^{-1}(1)$, i.e., $w_1$ satisfies
\beq
\lbl{16a}
\Delta w_1+kg_u(x,u)w_1=1, \s \mbox{in $D$}, \s w_1=0 \s \mbox{ on $\partial D$} \,.
\eeq
We also consider $z=1-L^{-1} \left( kg_u(x,u) \right)$, i.e., $z$ satisfies
\beq
\lbl{17}
\Delta z+kg_u(x,u)z=0, \s \mbox{in $D$}, \s z=1 \s \mbox{ on $\partial D$} \,.
\eeq
We shall build the solution of (\ref{16}), by using $w_1$ and $z$. Multiplying (\ref{16a}) by $z$, subtracting the equation (\ref{17}) multiplied by $w_1$, and integrating
\beq
\lbl{18}
\int _D z \, dx=\int_{\partial D} \frac{\partial w_1}{\partial n} \, ds \,.
\eeq

\noindent
{\bf Sub-case i}.  The integrals in (\ref{18}) are both non-zero. We now construct the solution of (\ref{16}) in the form
\[
w=\mu ^* (w_1+bz)+L^{-1} (\theta (x))+b_1z \,,
\]
with the constants $\mu ^*$, $b$ and $b_1$ to be selected. Clearly, $w$ satisfies the equation in (\ref{16}), while $w \,| \, _{\partial D}=b+b_1$. Since $\int _D z \, dx \ne 0$, we can select $b$ and $b_1$, such that
\[
\int _D (w_1+bz) \, dx=\int _D (L^{-1} \left( \theta )+b_1z \right) \, dx=0 \,.
\]
It follows that $\int _D w(x) \, dx=0$, for any $\mu ^*$.  The function $W \equiv w_1+bz$ satisfies
\beqa \nonumber
& \Delta W+kg_u(x,u)W=1, \s \mbox{in $D$} \\ \nonumber
& W \,| \, _{\partial D}=b,    \\ \nonumber
& \int _D W(x) \, dx=0 \,. \nonumber
\eeqa
By Lemma \ref{lma:2}, it follows that $\int_{\partial D} \frac{\partial W}{\partial n} \, ds \ne 0$. Hence, we can select $\mu ^*$, so that 
\[
\int_{\partial D} \frac{\partial w}{\partial n} \, ds=0 \,.
\]
\md 

\noindent
{\bf Sub-case ii}.  The integrals in (\ref{18}) are both zero. Since $\int _D z(x) \, dx=0$, it follows by Lemma \ref{lma:2} that
\beq
\lbl{19}
\int_{\partial D} \frac{\partial z}{\partial n} \, ds \ne 0 \,.
\eeq
Similarly, since $\int_{\partial D} \frac{\partial w_1}{\partial n} \, ds=0$, we conclude by Lemma \ref{lma:2} that
\beq
\lbl{20}
\int _D w_1(x) \, dx \ne 0 \,.
\eeq
We now construct the solution of (\ref{16}) in the form
\[
w=\mu ^* w_1+L^{-1} (\theta (x))+bz \,.
\]
Clearly, $w$ satisfies the equation in (\ref{16}), while $w \,| \, _{\partial D}=b$.
By (\ref{20}) we can choose $\mu ^*$, so that $\int _D w(x) \, dx=0$, for all $b$, and then we choose $b$, so that $\int_{\partial D} \frac{\partial w}{\partial n} \, ds=0$, by using (\ref{19}).
\bigskip 

\noindent
{\bf Case 2}. The operator $L$ is not invertible.
\md 

Since we have assumed $kg_u(x,u)<\la _2$, it follows that the null space of $L$ is one dimensional, spanned by some $\bar w $, with $\bar w (x)>0$ on $D$. By the Fredholm alternative, given any $f(x) \in L^2(D)$, $L^{-1}(f(x))$ is defined, i.e., the problem
\[
\Delta w+kg_u(x,u)w=f(x), \s \mbox{in $D$}, \s w=0 \s \mbox{ on $\partial D$}
\]
is solvable, if and only if $\int _D f(x) \bar w(x) \, dx=0$. We can choose a constant $s$, so that $\int _D \left(s+\theta(x) \right) \bar w(x) \, dx=0$, which implies that $L^{-1} \left(s+\theta(x) \right)$ is defined. Similarly, we select a constant $t$, so that  $L^{-1} \left(-kg_u(x,u) + t\right)$ is defined, and we set $z=1+L^{-1} \left(-kg_u(x,u) + t\right)$, i.e., $z$ satisfies
\beq
\lbl{21-}
\Delta z+kg_u(x,u)z=t, \s \mbox{in $D$}, \s z=1 \s \mbox{ on $\partial D$} \,.
\eeq

\noindent
{\bf Sub-case i}. $\int _D z \, dx \ne 0$. We  construct the solution of (\ref{16}) in the form
\[
w=L^{-1} \left(s+\theta(x) \right) +b_1 z+a (\bar w +b_2z) \,,
\]
with $s$  fixed above, and the constants $a$, $b_1$ and $b_2$ to be selected. We choose $b_1$ and $b_2$, so that
\[
\int _D \left[ L^{-1} \left(s+\theta(x) \right) +b_1 z \right] \, dx=\int _D (\bar w +b_2z) \, dx=0 \,.
\]
Then $\int _D w(x) \, dx=0$, for any constant $a$. The function $W(x)=\bar w +b_2z$ satisfies
\beqa \nonumber
& \Delta W+kg_u(x,u)W=b_2t, \s \mbox{in $D$} \\ \nonumber
& W \,| \, _{\partial D}=b_2,    \\ \nonumber
& \int _D W(x) \, dx=0 \,. \nonumber
\eeqa
By Lemma \ref{lma:2}, it follows that $\int_{\partial D} \frac{\partial W}{\partial n} \, ds \ne 0$. We select $a$ so that $\int_{\partial D} \frac{\partial w}{\partial n} \, ds = 0$. Then $w$ satisfies (\ref{16}). (Observe that $Lw=s+(b_1+a b_2)t+\theta (x) \equiv \mu ^*+\theta (x)$, $w \,| \, _{\partial D}=b_1+ab_2 \equiv b$.)
\md 

\noindent
{\bf Sub-case ii}. $\int _D z \, dx = 0$. By Lemma \ref{lma:2}, and (\ref{21-}), it follows that $\int_{\partial D} \frac{\partial z}{\partial n} \, ds \ne 0$. Let
\[
w(x)=L^{-1} \left(s+\theta(x) \right)+a \bar w  +b z\,.
\]
Choose the constant $a$, so that
\[
\int _D \left[ L^{-1} \left(s+\theta(x) \right)+a \bar w \right] \, dx=0  \,.
\]
Then, $\int _D w(x) \, dx=0$, for any $b$. By Lemma \ref{lma:2}, $\int_{\partial D} \frac{\partial z}{\partial n} \, ds \ne 0$. Then choosing $b$, so that  $\int_{\partial D} \frac{\partial w}{\partial n} \, ds = 0$, we conclude that $w(x)$  satisfies (\ref{16}). (Here  $Lw=s+bt+\theta (x) \equiv \mu ^*+\theta (x)$, $w \,| \, _{\partial D}= b$.)
\bigskip 

\noindent
{\bf Case II}: $\xi _1 \ne 0$. We reduce it to the case $\xi _1  =0$, by letting $v=u-\xi _1$. Then $v(x)$ satisfies

\beqa
\lbl{21}
& \Delta v+kg(x,v+\xi _1)=\mu+\theta (x) \s \mbox{in $D$} \\ \nonumber
& v \,| \, _{\partial D}=b- \xi _1, \s\s \int_{\partial D} \frac{\partial v}{\partial n} \, ds=0 \\ \nonumber
& \frac{1}{|D|} \int_{D} v(x) \, dx = 0 \,. \nonumber
\eeqa 
We perform the continuation in $k$ exactly the same way as before, since the bound on $|g_v(x,v+\xi _1)|$ remains the same.
\md 

We conclude that each solution of (\ref{10}) can be continued locally in $k$. To show that solutions can be continued for all $0 \leq k \leq 1$, we need an a priori bound on $(u,\mu)(k)$. Write $u(x)=\xi _1 +U(x)$, with $\int_{D} U(x) \, dx = 0$. We have
\beqa
\lbl{21.1}
& \Delta U+kg(x,\xi _1 +U)=p(x) \s \mbox{in $D$} \\
& U\,| \, _{\partial D}=b- \xi _1, \s\s \int_{\partial D} \frac{\partial U}{\partial n} \, ds=0 \,. \nonumber
\eeqa
By our conditions (\ref{11}) and (\ref{12}), we can find positive constants $c_1$ and $c_2$, with $c_1<c_0$, such that
\beq
\lbl{21.2}
|g(x,\xi _1 +U)| \leq c_1 U+c_2, \s \mbox{for all $x \in \bar D$, and $U \in R$} \,.
\eeq
Multiplying (\ref{21.1}) by $U$, and using the Lemma \ref{lma:0}, we conclude a bound on $||U||_{L^2(D)}$. Using this bound in (\ref{21.1}), and the estimate (\ref{21.2}), we obtain a bound on $||U||_{W^{2,2}(D)}$ (keep in mind that $\xi _1$ is arbitrary, but fixed). Writing $u=b+V$, and proceeding similarly, we conclude an $L^2$ bound on $V$, which implies a bound on $b$, since $||u||_{L^2(D)}$ is bounded.
\epf

Hence for any $\xi _1$ we have a curve of solutions $(u, \mu )(k) \in \bar H^2 \times R$ solving (\ref{10}), and at $k=1$, we have a solution of 
\beqa
\lbl{22}
& \Delta u+g(x,u)=\mu+\theta (x) \s \mbox{in $D$} \\ \nonumber
& u \,| \, _{\partial D}=b, \s\s \int_{\partial D} \frac{\partial u}{\partial n} \, ds=0 \\ \nonumber
& \frac{1}{|D|} \int_{D} u(x) \, dx = \xi _1 \,. \nonumber
\eeqa
Recall that the original problem (\ref{p1}), after we decomposed $p(x)=\mu _0+\theta(x)$, is
\beqa
\lbl{23}
& \Delta u+g(x,u)=\mu _0+\theta \s \mbox{in $D$} \\
& u \,| \, _{\partial D}=b, \s\s \int_{\partial D} \frac{\partial u}{\partial n} \, ds=0 \,. \nonumber
\eeqa
Hence, it remains to show that one can choose $\xi _1$, so that $\mu =\mu _0$ in (\ref{22}). The corresponding $u(x)$ (from (\ref{22})) is then a solution of (\ref{23}).
\md 

\noindent
{\bf Remark} If the domain $D$ is an interval in one dimension, then $c_0 =\la _2$. Indeed, suppose $D=(-L,L)$. Represent $u(x) \in W^{1,2}_0(-L,L)$ by its Fourier series
\[
u(x)=a_0+\Sigma _{n=1}^{\infty} a_n \cos \frac{n\pi}{L} x+b_n \sin \frac{n\pi}{L} x \,.
\]
If $\int _{-L}^L u(x) \, dx=0$, then $a_0 =0$, and hence
\[
\int _{-L}^L {u'}^2 (x) \, dx \geq \frac{\pi ^2}{L^2} \int _{-L}^L u^2(x) \, dx=0 \,,
\]
i.e., $c_0=\frac{\pi ^2}{L^2}$. The Dirichlet eigenvalues are $\la _n=\frac{n^2\pi ^2}{(2L)^2}$, so that $c_0 =\la _2$.

\section{Continuation in $\xi _1$}
\setcounter{equation}{0}

We have just seen that for each $\xi _1$ 
the problem
\beqa
\lbl{24}
& \Delta u+g(x,u)=\mu+\theta (x) \s \mbox{in $D$} \\ \nonumber
& u \,| \, _{\partial D}=b, \s\s \int_{\partial D} \frac{\partial u}{\partial n} \, ds=0
\eeqa
has a solution $(u, \mu) \in \bar H^2 \times R$, with the average value of $u(x)$ equal to $\xi _1$. We now show that all these solutions lie on a unique solution curve, which is globally parameterized by  $\xi _1$, the average value of the solution ($\xi _1=\frac{1}{|D|} \int_{D} u(x) \, dx$).

\begin{thm}\lbl{thm:2} 
Assume that the conditions of the Theorem \ref{thm:1} hold. Then any solution of (\ref{24}) can be continued in $\xi _1$, for all $-\infty<\xi _1<\infty$, giving us a curve of solutions $(u,b,\mu)(\xi _1) \in  W^{2,2}(D) \times R \times R$. Moreover, for each  $\xi _1$ there exists a unique solution $(u,b,\mu)$ of (\ref{24}). All solutions of the problem (\ref{24}) lie on a unique continuous solution curve $(u,b,\mu)(\xi _1) \in W^{2,2}(D) \times R \times R$, with $\xi _1$ being a global parameter.
\end{thm}

\pf
We use the Implicit Function Theorem to show that any  solution of (\ref{24}) can be continued in $\xi _1$. The proof is essentially the same, as the one above for continuation in $k$. Defining the map
\[
F(v,\mu ,\xi _1)=\Delta v + g(x,v+\xi _1)-\mu \; : \; \bar H^2 \times R  \times R \ra L^2(D) \,,
\]
we recast  the problem (\ref{24})  in the operator form
\beq
\lbl{25}
F(v,\mu ,\xi _1)=\theta (x) \,.
\eeq
Similarly to the Theorem \ref{thm:1}, we apply the  Implicit Function Theorem to perform continuation in $\xi _1$, by  showing  that the linearized operator
\[
F_{(v,\mu)}(u,\mu ,\xi _1)(w, \mu ^*)=\Delta w+g_u(x,v+\xi _1)w-\mu ^*
\]
is one-to-one and onto.
\medskip

Assume that there are two different pairs of solutions $(u_1(x), \mu_1)$ and $(u_2(x), \mu_2)$, with $u_1(x)$ and $u_2(x)$ having the same average value $\xi _1^0$. These pairs are also solutions of (\ref{10}), with $k=1$. By the Theorem \ref{thm:1}, we can continue both pairs for decreasing $k$ on two curves of solutions, with fixed average $\xi _1^0$. These curves do not intersect, since at all points the Implicit Function Theorem applies. At $k=0$, we obtain a contradiction with the uniqueness part of Lemma \ref{lma:2}.
\epf

\section{Existence  of solutions}
\setcounter{equation}{0}

We now discuss which values of $\mu$ are covered, as we continue in $\xi _1 \in (-\infty,\infty)$,  which translates into existence results for our problem (\ref{p1}).

\begin{thm}\lbl{thm:3}
Assume that the function $g(x,u) \in L^{\infty}(\bar D \times R)$ satisfies the conditions (\ref{p4}), (\ref{11}) and (\ref{12}). Then the condition (\ref{p5}) is necessary and sufficient for the existence of solution for the problem (\ref{p1}).
\end{thm}

\pf
We proved above that  the condition (\ref{p5}) is necessary for solvability, we now prove its sufficiency. By the Theorem \ref{thm:1}, for any $\xi _1 \in R$, there exists a triple $(u,b,\mu)=(u,b,\mu)(\xi _1)$ solving the problem (\ref{p1}), with average of $u(x)$ equal to $\xi _1$. As before, we write $u=\xi _1+U$, with $U=U(\xi _1)$ satisfying $\int _D U \, dx=0$. By the Theorem \ref{thm:2},
$\mu=\mu(\xi _1)$ is continuous in $\xi _1$, and integrating the equation (\ref{p1}), we express
\beq
\lbl{26}
\mu =\frac{1}{|D|} \int_{D} g \left(x,\xi _1+U(\xi _1) \right) \, dx \,.
\eeq
 We need to show that we can choose $\xi _1=\xi^0 _1$, so that $\mu(\xi^0 _1)=\mu _0$. This will follow, once we prove the existence of the limits
\beq
\lbl{27}
\lim _{\xi _1 \ra \infty} \int_{D} g \left(x,\xi _1+U(\xi _1) \right) \, dx =\int_{D} g \left(x,\infty \right) \, dx \,,
\eeq
\beq
\lbl{28}
\lim _{\xi _1 \ra -\infty} \int_{D} g \left(x,\xi _1+U(\xi _1) \right) \, dx =\int_{D} g \left(x,-\infty \right) \, dx \,.
\eeq
In view of (\ref{26}), $U=U(\xi _1)$ is satisfying
\beqa \nonumber
& \Delta U=-g \left(x,\xi _1+U \right)+\frac{1}{|D|} \int_{D} g \left(x,\xi _1+U \right) \, dx+\theta (x), \s \mbox{in $D$} \\ \nonumber
& U \,| \, _{\partial D}=b_1, \s\s \int_{\partial D} \frac{\partial U}{\partial n} \, ds=0   \\ \nonumber
&  \frac{1}{|D|} \int _D U(x) \, dx=0 \,, \nonumber
\eeqa
with $b_1=b-\xi _1$, a new unknown constant. By Corollary 1 to Lemma \ref{lma:1}, there is a constant $c$, so that
\[
||U(\xi _1)||_{L^{\infty} (\bar D)} \leq c, \s \mbox{uniformly in $\xi _1 \in R$} \,,
\]
from which the limits (\ref{27}) and (\ref{28}) follow.
\epf

The following result is similar to the one in D.G. de Figueiredo and W.-M.  Ni \cite{FN}, who considered the case of Dirichlet problem.

\begin{thm}\lbl{thm:4}
Consider the problem
\beqa 
\lbl{29}
& \Delta u+g(x,u)=\theta(x) \s \mbox{in $D$} \\ \nonumber
& u \,| \, _{\partial D}=b, \s\s \int_{\partial D} \frac{\partial u}{\partial n} \, ds=0 \,, \nonumber
\eeqa
where $\theta(x) \in L^2(D)$ is a given function, satisfying $\int _D \theta(x) \, dx=0$.
Assume that the function $g(x,u) \in L^{\infty}(\bar D \times R)$ satisfies the conditions (\ref{11})  and (\ref{12}), and in addition
\beq
\lbl{30}
ug(x,u)>0, \s\s \mbox{for all $u \in R$, and $x \in D$} \,.
\eeq
Then the problem (\ref{29}) has a solution $(u,b) \in W^{2,2}(D) \times R$.
\end{thm}

\pf
The proof is similar to the one above. Here we have $\mu ^0 =0$. Using (\ref{30}), we see from (\ref{26}) that $\mu >0$ ($<0$) for $\xi _1$ large and positive (negative).
By continuity, we have $\mu (\xi ^0_1)=0$ for some $\xi ^0_1$.
\epf

We also have the following multiplicity result.
\begin{thm}\lbl{thm:5}
Consider the problem
\beqa 
\lbl{31}
& \Delta u+g(x,u)=\mu+\theta(x) \s \mbox{in $D$} \\ \nonumber
& u \,| \, _{\partial D}=b, \s\s \int_{\partial D} \frac{\partial u}{\partial n} \, ds=0 \,, \nonumber
\eeqa
where $\theta(x) \in L^2(D)$ is a given function, satisfying $\int _D \theta(x) \, dx=0$. Assume there is a constant $G>0$ such that 
\beq
\lbl{32}
ug(x,u)>0, \s\s \mbox{for all $|u|>G$, and $x \in D$} \,.
\eeq
Assume that the limits $ g(x, \pm \infty)=\lim _{u \ra \pm \infty} g(x,u) $ exist, uniformly in $x \in \bar D$, and 
\beq
\lbl{33}
g(x, \pm \infty) \equiv 0 \,.
\eeq
Then there exists an interval $(\mu _-,\mu _+)$, with $\mu _- <0$ and $\mu _+ >0$, so that for $\mu \in (\mu _-,\mu _+) \setminus \{0\}$ the problem (\ref{31}) has at least two solutions $(u_1,b_1)$ and $(u_2,b_2)$.
\end{thm}

\pf
By the Theorem \ref{thm:1}, for any $\xi _1 \in R$, there exists a triple $(u,b,\mu)=(u,b,\mu)(\xi _1)$ solving the problem (\ref{31}), with average of $u(x)$ equal to $\xi _1$. As in Theorem \ref{thm:3}, we continue this solution in $\xi _1$, paying particular attention to the curve $\mu=\mu(\xi _1)$.
From the formula (\ref{26}) and the condition (\ref{32}) it follows that $\mu(\xi _1)$ is positive (negative) for $|\xi _1|$ large and positive (negative). 
The condition (\ref{33}) implies that $\lim _{|\xi _1| \ra \infty}\mu(\xi _1)=0$, i.e., both sides of the curve $\mu=\mu(\xi _1)$ tend to zero.
\epf

\end{document}